\documentclass[reqno,A4paper]{amsart}

\usepackage{amsmath}
\usepackage{amssymb}
\usepackage{amsthm}
\usepackage{enumerate}
\usepackage{graphicx}
\usepackage{mathrsfs} 
\usepackage{eqlist}
\usepackage{array}

\theoremstyle{plain}
\newtheorem{theorem}{Theorem}[section]

\newtheorem{lemma}{Lemma}[section]
\newtheorem{corollary}{Corollary}[theorem]

\theoremstyle{definition}
\newtheorem{definition}{Definition}
\newtheorem{remark}{\textup{Remark}} 

\newtheorem{question}{\textit{Question}}

\numberwithin{equation}{section}

\begin{document}

\title[Properties of Certain Analytic Functions]%
{Theory of certain Non-Univalent Analytic functions}
\author[Kamaljeet Gangania]%
{Kamaljeet Gangania}

\newcommand{\acr}{\newline\indent}

\address{Department of Applied Mathematics\acr
                   Delhi Technological University\acr
                   New Delhi\acr
                   India}
\email{gangania.m1991@gmail.com}


\thanks{The work of Kamajeet Gangania is supported by University Grant Commission, New-Delhi, India  under UGC-Ref. No.:1051/(CSIR-UGC NET JUNE 2017).}


\begin{abstract}
We investigate the non-univalent function's properties reminiscent of the theory of univalent starlike functions.
Let the analytic function $\psi(z)=\sum_{i=1}^{\infty}A_i z^i$, $A_1\neq0$ be univalent in the unit disk. Non-univalent functions may be found in the class $\mathcal{F}(\psi)$ of analytic functions $f$ of the form $f(z)=z+\sum_{k=2}^{\infty}a_k z^k$ satisfying  $({zf'(z)}/{f(z)}-1) \prec \psi(z)$. Such functions, like the Ma and Minda classes of starlike functions, also have nice geometric properties. For these functions, growth and distortion theorems have been established. Further, we obtain bounds for some sharp coefficient functionals and establish the Bohr and Rogosinki phenomenon for the class $\mathcal{F}(\psi)$. Non-analytic functions that share properties of analytic functions are known as Poly-analytic functions. Moreover, we compute Bohr and Rogosinski's radius for Poly-analytic functions with analytic counterparts in the class $\mathcal{F}(\psi)$ or classes of Ma-Minda starlike and convex functions. 
\end{abstract}

\maketitle

	\section{Introduction}
\label{intro}

Let $\mathcal{A}$ be the class of analytic functions of the form $f(z)=z+\sum_{k=2}^{\infty}a_kz^k$ in the open unit disk ${\Delta}:=\{z: |z|<1\}$.
Let $f(z)=w$ and $\Gamma_w$ be the image of $\Gamma_z$ (the circle $C_r: z=re^{i\theta}$) under the function $f$ in $\mathcal{A}$. The curve $\Gamma_w$ is said to be starlike with respect to $w_0=0$ if $\arg(w-w_0)$ is a non-decreasing function of $\theta$, that is,
\begin{equation*}\label{arg-def}
\frac{d}{d\theta} \arg(w-w_0)\geq0, \quad \theta\in [0,2\pi],
\end{equation*}  
which is equivalent to
\begin{equation}\label{charcter}
\frac{d}{d\theta} \arg(w-w_0)=\Re\left(\frac{zf'(z)}{f(z)}\right)\geq0.
\end{equation}

From \eqref{charcter}, we note the importance of the Carath\'{e}odory functions when \eqref{charcter} is viewed in terms of subordination as:
\begin{equation}\label{star-subord}
\frac{zf'(z)}{f(z)}\prec \frac{1+z}{1-z} \quad (z\in\Delta),
\end{equation}
where the symbol $\prec$ means for the usual subordination. In view of \eqref{star-subord}, let $\Psi$ be an analytic function in $\Delta$ with positive real part, $\Psi({\Delta})$ symmetric about the real axis with $\Psi'(0)>0$ and $\Psi(0)=1$.  Generalizing \eqref{star-subord}, Ma and Minda~\cite{minda94} considered the unified approach by defining the class of starlike functions
\begin{equation}\label{mindaclass}
\mathcal{S}^*(\Psi):= \biggl\{f\in \mathcal{A} : \frac{zf'(z)}{f(z)} \prec \Psi(z) \biggl\}.
\end{equation}

If the inequality \eqref{charcter} holds for each circle $|z|=r<1$, then it characterizes a special class $\mathcal{S}^{*}$, the class of starlike functions in ${\Delta}$. It is obvious from \eqref{arg-def} that for each $0<r<1$, the curve $\Gamma_w$ is not allowed to have a loop which ensure the univalency of the function. But if for some $0\neq z\in \Delta$, $\Re(zf'(z)/f(z))<0$, then $f$ is not starlike with respect to $0$, or equivalently we can say that the image curve $\Gamma_w: f(|z|=r)$ is definitely not starlike, but still it may or may not be univalent.

Therefore, view of \eqref{charcter}, \eqref{star-subord} and the class \eqref{mindaclass}, if we take
\begin{equation} \label{new-character}
\frac{zg'(z)}{g(z)} -1\prec \psi(z),
\end{equation}
{\it{where $\psi$ be the analytic univalent function in ${\Delta}$ such that $\psi(0)=0$, $\psi({\Delta})$ is starlike with respect to $0$}}. Then $g$ satisfying \eqref{new-character} may not be univalent. In view of \eqref{new-character}, Kargar et al.~\cite{kargar-2019} studied a class reminiscent to the Ma-Minda class.  
\begin{equation*}\label{boothlem}
\mathcal{BS}(\beta):= \biggl\{f\in \mathcal{A} : \frac{zf'(z)}{f(z)}-1 \prec \frac{z}{1-\beta z^2},\; \beta\in [0,1) \biggl\},
\end{equation*}
where  $z/(1-\beta z^2)=:\psi(z)$ (Booth lemniscate function~\cite{piejko-2013} and \cite{piejko-2015}) is an analytic univalent function and symmetric with respect to the real and imaginary axes.  Also, see the family of generalized pascal snails~\cite{Kanas-pascal}. Masih et al. \cite{Masih-2019} considered the following class with $\beta\in [0,1/2]$:
\begin{equation*}\label{cissoidclass}
\mathcal{S}_{cs}(\beta):= \biggl\{f\in \mathcal{A} : \left(\frac{zf'(z)}{f(z)}-1\right) \prec \frac{z}{(1-z)(1+\beta z)},\; \beta\in [0,1) \biggl\}.
\end{equation*}
Kumar and Gangania~\cite{ganga-geo2021} studied the class
\begin{equation*}
\mathcal{S}_{\gamma}(\eta):= \biggl\{f\in \mathcal{A} : \left(\frac{zf'(z)}{f(z)}-1\right) \prec \frac{\gamma z}{(1+\eta z)^2},\; \eta\in [0,1),\; \gamma>0 \biggl\}. 
\end{equation*}

Motivated by the above classes and observations, Kumar and Gangania~\cite{ganga-geo2021} introduced the following class:
\begin{definition}\cite{ganga-geo2021}\label{FPsi-2}
	Let  $\psi$ be the analytic univalent function in ${\Delta}$ such that $\psi(0)=0$, $\psi({\Delta})$ is starlike with respect to $0$. Then
	\begin{equation*}\label{gen-ma-min}
	\mathcal{F}(\psi):= \left\{f\in \mathcal{A}:  \left( \frac{zf'(z)}{f(z)}-1 \right) \prec \psi(z),\; \psi(0)=0 \right\}.
	\end{equation*}
\end{definition}

The following function plays the role of extremal
\begin{equation}\label{f0}
f_0(z)=z \exp\int_{0}^{z}\frac{\psi(t)}{t}dt.
\end{equation}

Notice that for $1+\psi(z)\not \prec (1+z)/(1-z)$, the functions in the class $\mathcal{F}(\psi)$ may not be univalent in ${\Delta}$ which also implies $\mathcal{F}(\psi)\not\subseteq \mathcal{S}^{*}$. Authors in \cite{ganga-geo2021} noticed that such non-univalent analytic functions also exhibit nice geometric properties similar to class $\mathcal{S}^*(\Psi)$ such as
\begin{theorem}\cite{ganga-geo2021}(Covering theorem)
	Let $\min_{|z|=r}\Re\psi(z)=\psi(-r)$. If  $f\in \mathcal{F}(\psi)$ and $f_0$ as defined in \eqref{f0}, then either $f$ is a rotation of $f_0$ or
	$$	\{w\in \mathbb{C} : |w|\leq-{f}_0(-1) \} \subset f({\Delta}),$$
	where $-{f}_0(-1)=\lim_{r\rightarrow 1}(-f_0(-r)).$
\end{theorem}
\begin{remark}
	{\it The class $\mathcal{F}(\psi)$ gives a new insight and a pool of opportunities to study non-univalent functions. Several problems, namely Area estimates, Logarithmic coefficient bounds, sharp coefficient bounds, hankel determinant, Zalcman functional, and several radius problems like majorization, radius of convexity, Bohr-phenomenon, differential subordination implication, preservation under convolution and some famous integral etc.} are yet to be known. 
\end{remark}

\begin{figure}[h]
	\begin{tabular}{c}
		\includegraphics[scale=0.5]{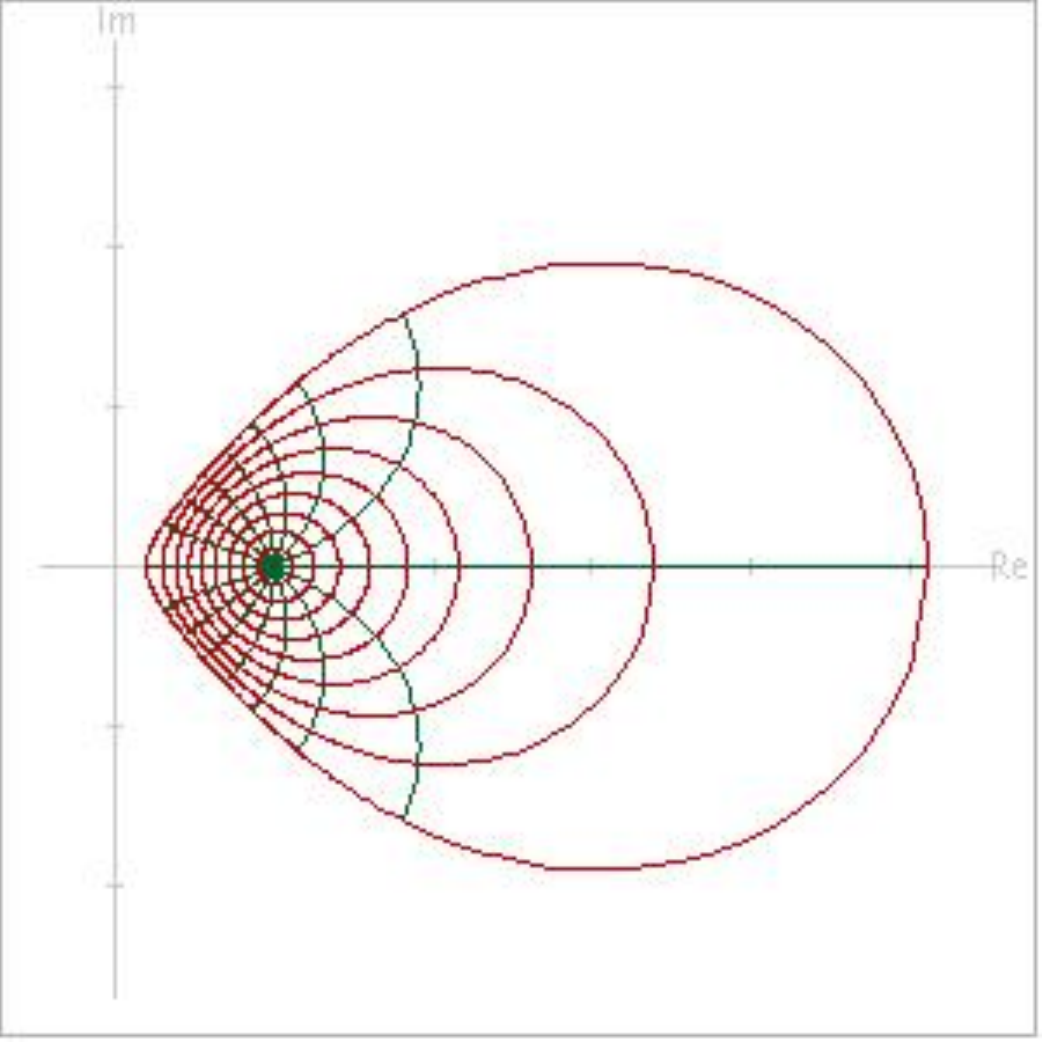}
	\end{tabular}
	\begin{tabular}{l}
		\includegraphics[scale=0.5]{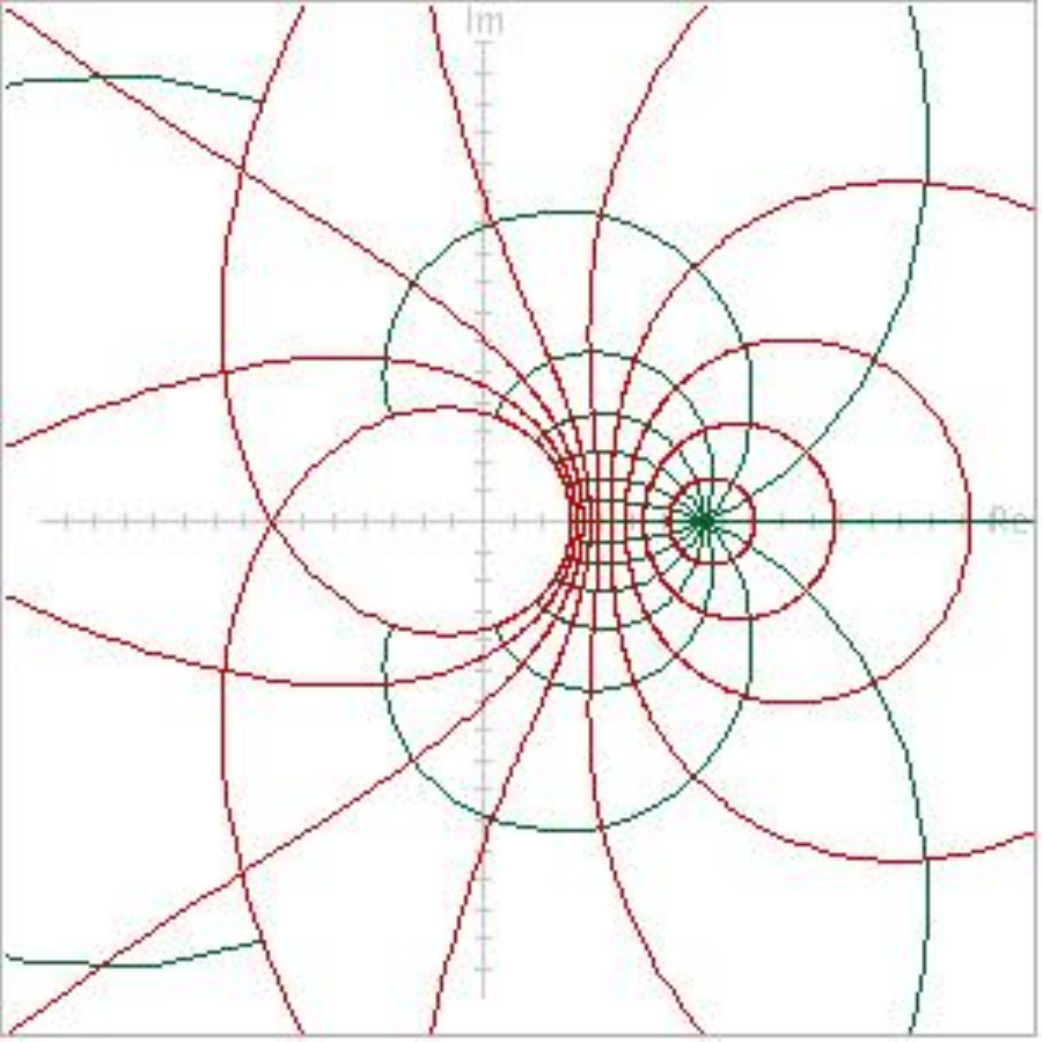}
	\end{tabular}
	\caption{Left side graph is for $\hat{f}(z)/z$ $(\beta=0.9)$ and Right side graph is for $\kappa(z)/z$ $(\eta=0.9)$.} \label{noncaratheodorygraph}
\end{figure}

Moreover, inclusion relation like $\mathcal{S}^*(\Psi) \subseteq \mathcal{S}^*$ is not available for the class $\mathcal{F}(\psi)$. This fact makes the above discussed problems non-trivial in $\mathcal{F}(\psi)$. In fact, several implications which holds for $\mathcal{S}^*(\Psi)$ are not true for the class $\mathcal{F}(\psi)$. For instance, if $f\in \mathcal{S}^*(\Psi)$ then $\Re(f(z)/z)>0$ for all $z\in \Delta$, see~\cite{minda94}. But $\Re(f(z)/z)>0$ may not be true for $f\in \mathcal{F}(\psi)$. See Figure~\ref{noncaratheodorygraph}. The function 
$$\hat{f}(z)=z\left(\frac{1+z\sqrt{\beta}}{1-z\sqrt{\beta}}\right)^{\frac{1}{2\sqrt{\beta}}} \in \mathcal{BS}(\beta) \Rightarrow \Re \left(\frac{\hat{f}(z)}{z} \right)>0\quad \text{for all}\; z\in \Delta.$$ 
However, for all $z\in \Delta$, the function
$${\kappa}(z):=z \exp\left(\frac{z}{(1+\eta z)^2}\right) \in \mathcal{S}_{1}(\eta) \Rightarrow \Re\left(\frac{\kappa(z)}{z} \right)\not>0, \quad \eta >2-\sqrt{3}.$$

Likewise, Non-analytic functions with significant structure and properties similar to those satisfied by analytic functions exist as well. Non-analytic functions with such a nice structure are referred to as poly-analytic functions, which we shall discuss briefly soon in Section~\ref{Radius-problems}. For there interesting historical viewpoint and the application's side, we refer to read \cite[Page~170]{Abdulhadi-Hajj2022}

The version of Landau's theorem for poly-analytic functions is given in
\begin{theorem}\cite[Theorem~1]{Abdulhadi-Hajj2022}(Landau's theorem)\label{1}
	Let $F(z)=\sum_{k=0}^{\alpha-1}{\bar{z}}^k f_{k}(z)$ be a poly-analytic function of order $\alpha$, $\alpha \geq 2$ on $\Delta$, where $f_k$ are analytic such that $f_k(0)=0$, $f'_{k}(0)=1$ and $|f_k(z)| \leq M$, for all $k$, with $M>1$. Then there is a constant $0< \rho_1 <1$ so that $F$ is univalent in $|z|< \rho_1.$ In particular, $\rho_1$ satisfies
	\begin{equation*}
	1-M \left( \frac{\rho_1(2-\rho_1)}{(1-\rho_1)^2} + \sum_{k=1}^{\alpha-1} \frac{\rho_{1}^{k}(1+k+k\rho_1)}{(1-k\rho_1)^2} \right)=0
	\end{equation*}
	and $F(\Delta_{\rho_1})$ contains a disk $\Delta_{R_1}$, where
	\begin{equation*}
	R_1= \rho_1-{\rho_1}^2 \left( \frac{1-{\rho_1}^{\alpha-1}}{1-\rho_1} \right) -M \sum_{k=0}^{\alpha-1} \frac{{\rho_1}^{k+2}}{1-\rho_1}.
	\end{equation*}
\end{theorem}  
For other relevant work in this direction, we refer to see \cite{Abdulhadi-Hajj2022,ChenGH-Bloch2000,Liu-Bloch2009} and references therein.

For poly-analytic functions of order $\alpha$ given by $F(z)=\sum_{k=0}^{\alpha-1} {\bar{z}}^k f_k(z)$, where $f_k(z)= \sum_{n=0}^{\infty} a_{n,k} z^n$, Abdulhadi and El Hajj~\cite{Abdulhadi-Hajj2022} defined the majorant series $M(F, r)$ by
$$ M(F, r)= \sum_{n=0}^{\infty}\sum_{k=0}^{\alpha-1}|a_{n,k}|r^{k+n}, $$ 
and proved the following Bohr's type inequality
\begin{theorem}\cite[Theorem~2]{Abdulhadi-Hajj2022}\label{2}
	Let $F(z)=\sum_{k=0}^{\alpha-1}{\bar{z}}^k f_k(z)$ be a poly-analytic function of order $\alpha$, where $f_k$ are analytic mappings for $k=0, 1\cdots, \alpha-1$, such that $g_k=f_0+\bar{f_k}$ is sense-preserving in $\Delta$ and $f_k(0)=0$ for each $k=0, \cdots, \alpha-1$. Suppose that $f_0$ is univalent and normalized by $f_0(0)=0$, $f'_{0}(0)=1$. Then 
	$$M(F,r) <1 \quad if \quad |z|<r_0,$$
	where $r_0\approx 0.318$ is the minimum positive root of the polynomial	
	$$(1-r)^3-r+r^{\alpha+1}=0.$$
\end{theorem}

The classical Bohr's inequality~\cite{bohr1914} says that if $f(z)=\sum_{n=0}^{\infty}a_n z^n$ is analytic in $\Delta$ and $|f(z)|<1$ for all $z\in \Delta$, then $\sum_{n=0}^{\infty}|a_n z^n|<1$ for all $z$ in $\Delta $ with $|z|<1/3$. The constant $1/3$ can not be improved and is known as Bohr's radius. In analogy with Bohr's inequality, there is also the notion of Rogosinski radius, however a little is known about Rogosinski radius as compared to Bohr radius, which states that, see \cite{Landau-1986,Rogosinski-1923,Schur-1925}: if $g(z)=\sum_{k=0}^{\infty}b_k$ with $|f(z)|<1$, then for every $N\geq1$ we have
$| \sum_{k=0}^{N-1}b_k z^k| \leq1$ for all $|z|\leq {1}/{2}$. The radius $1/2$ is called the Rogosinski radius.

Kayumov et al.~\cite{Kayumov-2021} studied a new quantity, called Bohr-Rogosinski sum and is given below:
\begin{equation*}
|g(z)|+ \sum_{k=N}^{\infty}|b_k||z|^k, \quad |z|=r.
\end{equation*}
For the case $N=1$, note that this sum is similar to the Bohr's sum, where $g(0)$ is replaced by $|g(z)|$. We also refer the readers to see~\cite{Aizenberg-2012,Alkha-2020}. Now we say the family $S(f)$ has Bohr-Rogosinski phenomenon, if there exists $r^{f}_{N} \in (0,1]$ such that the inequality: 
$$|g(z)|+ \sum_{k=N}^{\infty}|b_k||z|^k \leq |f(0)|+ d(f(0),\partial \Omega)$$
holds for $|z|=r \leq r^{f}_{N}$, where $d(f(0),\partial\Omega)$ denotes the Euclidean distance between $f(0)$ and the boundary of $\Omega=f(\Delta)$. The largest such $r^{f}_{N}$ is called the Bohr-Rogosinski radius.

In this investigation, we obtain sharp initial coefficient's estimates along with some well known coefficient functionals for the class $\mathcal{F}(\psi)$. We also derive Bohr and Bohr-Rogosinski radius for $\mathcal{F}(\psi)$. Further, we explore the classes $\mathcal{F}(\psi)$ and $\mathcal{S}^*(\Psi)$ in connection with the poly-analytic functions to obtain Bohr and Rogosinski's type radius.
\section{Radius problems: The Bohr and Rogosinski phenomenon}\label{Radius-problems}
For our next result, we need to recall the following result of  Ruscheweyh and Stankiewicz~\cite{RusStan-1985}:
\begin{lemma}[\cite{RusStan-1985}]\label{convo-result}
	Let the analytic functions $F$ and $G$ be convex univalent in ${\Delta}$. If $f\prec F$ and $g \prec G$, then $$f*g \prec F*G \quad (z\in{\Delta}).$$
\end{lemma}

We also refer to see~\cite{minda94} for a parallel result concerning the class $\mathcal{S}^*(\Psi)$.
\begin{lemma}\label{subordination}
	Let $\psi(\Delta)$ be a convex domain . If $f \in \mathcal{F}(\psi)$. Then $f(z)/z \prec f_{\psi}(z)/z$, where
	$f_{\psi}(z)=z\exp \int_{0}^{z}(\psi(t)/t)dt$.
\end{lemma}
\begin{proof}
	Let $f\in \mathcal{F}(\psi)$. Then
	\begin{equation}\label{phi-psi}
	\phi(z):= \frac{zf'(z)}{f(z)}-1 \prec \psi(z).
	\end{equation}
	Let us take the univalent convex function 
	\begin{equation*}
	g(z)=\log \left(\frac{1}{1-z}\right) =\sum_{n=1}^{\infty}\frac{z^n}{n}.
	\end{equation*}
	Thus for  $f\in \mathcal{A}$, we get
	\begin{equation}\label{phi-gExp1}
	\phi(z)* g(z)= \int_{0}^{z}\frac{\phi(t)}{t}dt \quad \text{and}\quad \psi(z)* g(z)= \int_{0}^{z}\frac{\psi(t)}{t}dt. 
	\end{equation}
	Since $\psi$ is convex. Therefore, applying Lemma~\ref{convo-result} in \eqref{phi-psi}, we get
	\begin{equation}\label{phi-gExp2}
	\phi(z)* g(z) \prec \psi(z)*g(z).
	\end{equation}
	Now from \eqref{phi-gExp1} and \eqref{phi-gExp2}, we obtain
	\begin{equation*}
	\int_{0}^{z}\frac{\phi(t)}{t}dt \prec \int_{0}^{z}\frac{\psi(t)}{t}dt, 
	\end{equation*}
	which  implies that
	\begin{equation*}
	\frac{f(z)}{z}:= \exp\int_{0}^{z}\frac{\phi(t)}{t}dt \prec \exp\int_{0}^{z}\frac{\psi(t)}{t}dt=: \frac{f_{\psi}(z)}{z}.
	\end{equation*}
	This completes the proof. 	
\end{proof}	

\begin{remark}{(Behaviour of the functions $f(z)/z$)}
	It is worth to mention that the extremal function $f_{\psi}(z)/z$ may not be univalent, Carath\'{e}odory and convex in Lemma~\ref{subordination}, see~\cite{ganga-geo2021} and Figure~\ref{noncaratheodorygraph}. This restrict us to examine the  parallel results as in the Ma and Minda classes. But still this information is useful to settle the Bohr-phenomenon for the class $\mathcal{F}(\psi)$. 	
\end{remark}

To discuss the Bohr phenomenon for the class $\mathcal{F}(\psi)$, we recall the following which is the analogue of the growth theorem~\cite{minda94} for $\mathcal{S}^*(\Psi)$:
\begin{lemma}\cite{ganga-geo2021}(Growth Theorem)\label{gen-thm1}
	Let us suppose that $\max_{|z|=r}\Re\psi(z)=\psi(r)$ and $\min_{|z|=r}\Re\psi(z)=\psi(-r)$. Then $f\in \mathcal{F}(\psi)$ satisfies the sharp inequalities:
	\begin{equation}\label{maingththm-eq}
	r \exp\left(\int_{0}^{r}\frac{\psi(-t)}{t}dt\right) \leq |f(z)| \leq
	r \exp\left(\int_{0}^{r}\frac{\psi(t)}{t}dt\right), \quad (|z|=r).
	\end{equation}
\end{lemma}

Without loss of generality, we can have in Lemma~\ref{gen-thm1} that $\max_{|z|=r}\Re\psi(z)=M(r)$ and $\min_{|z|=r}\Re\psi(z)=m(r)$, where $M$ and $N$ are functions of $r$ for each $z$ in $\Delta$ such that $|z|=r$.
\begin{remark}\label{generaluppergrowth}
	In case of carath\'{e}odory functions $p$ and $q$, we see that $p\prec q$, where $q(z)=(1+z)/(1-z)$ is the best dominant. But in case of functions in $\mathcal{S}$, we have no such best dominant. Therefore, it is worthy to obtain upper bound for $|f(z)|$ in general for the class $\mathcal{F}(\psi)$, where $\psi\in \mathcal{S}$. Note that
	\begin{equation*}
	\max_{|z|=r} \left|\frac{zf'(z)}{f(z)} \right| \leq \frac{1-r+r^2}{(1-r)^2}.
	\end{equation*}
	This yields
	\begin{equation*}
	|f'(z)| \leq \frac{(1-r+r^2)}{(1-r)^2} \exp \left(\frac{r}{1-r} \right) \quad \text{and} \quad |f(z)| \leq r \exp\left(\frac{r}{1-r}\right).
	\end{equation*}
	The equality cases hold for the non-univalent function $$f(z)= z\exp \left(\frac{z}{1-z} \right).$$ 
\end{remark}

At this conjunction in view of the remark~\ref{generaluppergrowth}, we propose the following:
\begin{question}
	Let $m(r):= \min_{|z|=r}|f(z)|$, where $f\in \mathcal{F}(\psi)$ and $\psi\in \mathcal{S}$. Find the value of $m(r)$.
\end{question}	

In the following result, it is important to note that bounds for $|b_n|$ are not available.	In rest of this section, we assume that coefficient $A_i$ of $\psi$ are positive.
\begin{remark}
In case, when all the coefficients $A_i$ of function $\psi$ are not positive, then without loss of generality we can appropriately replace the function $f_\psi(z)=z \exp \int_{0}^{z}(\psi(t)/t) dt:=z+\sum_{n=2}^{\infty}a_n z^n$ by the analytic function $\hat{f}_{\psi}(z):= z+\sum_{n=2}^{\infty}|a_n|z^n$, throughout this paper.
\end{remark}
\begin{theorem}
	Let $\psi(\Delta)$ be convex. The class $S(f)$, where  $f \in \mathcal{F}(\psi)$ satisfies the Bohr-phenomenon
	\begin{equation*}
	\sum_{n=1}^{\infty}|b_n| |z|^n \leq d(0, \partial f(\Delta))
	\end{equation*}
	in $|z|=r\leq \min\{1/3, r_0 \}$, where $r_0$ is least positive root of the equation
	\begin{equation}\label{boothbohr-eq}
	f_{\psi}(r)- \exp \int_{0}^{1}\frac{\psi(-t)}{t}dt=0.
	\end{equation}
	The result is sharp in case of $r_0 \leq 1/3$.
\end{theorem}
\begin{proof}
	Since  $g\in S(\mathcal{F}(\psi))$, we have $g\prec f$ for a fixed $f\in \mathcal{F}(\psi)$. From Lemma~\ref{gen-thm1}, we obtain the Koebe-radius 
	$$r_{*}=\exp \int_{0}^{1}\frac{\psi(-t)}{t}dt$$ 
	such that $  r_{*}\leq d(0,\partial\Omega)=|f(z)|$ for $|z|=1$. Also using Lemma~\ref{subordination}, we have
	\begin{equation}\label{f-f0}
	\frac{f(z)}{z}\prec \frac{f_\psi(z)}{z}.
	\end{equation}
	Recall the result \cite[Lemma~2.1]{kamal-mediter2021}, which reads as:
	let $f$ and $g$ be analytic in ${\Delta}$ with $g\prec f,$ where
	$f(z)=\sum_{n=0}^{\infty}a_n z^n$ and $ g(z)= \sum_{k=0}^{\infty}b_k z^k.$
	Then
	$\sum_{k=0}^{\infty}|b_k|r^k \leq \sum_{n=0}^{\infty}|a_n|r^n$ for $ |z|=r\leq1/3.$
	Now using this result for  $g\prec f$ and  \eqref{f-f0}, we have
	\begin{equation*}
	\sum_{k=1}^{\infty}|b_k|r^k \leq	r+\sum_{n=2}^{\infty}|a_n|r^n \leq f_{\psi}(r)\quad\text{for}\; |z|=r\leq1/3.
	\end{equation*}
	Finally, to establish the inequality
	$\sum_{k=1}^{\infty}|b_k|r^k \leq d(f(0),\partial\Omega),$
	it is enough to show $f_{\psi}(r) \leq r_{*}.$
	But this holds whenever $r\leq r_0$, where $r_0$ is the least positive root of the equation $f_{\psi}(r)=r_{*}.$ Now let 
	$$G(r):=f_{\psi}(r)-r_{*}.$$
	Note that $G'(r)=f'_{\psi}(r)>0$ whenever the Taylor coefficient's of $f_{\psi}$ are positive, yields unique positive root of $G(r)=0$. Also, $G$ is a continuous function of $r$ in $(0,1)$ such that
	\begin{equation*}
	G(0)=f_{\psi}(0)-\exp\int_{0}^{1}\frac{\psi(-1)}{t}dt <0
	\end{equation*}
	and 
	\begin{equation*}
	G(1)=f_{\psi}(1)-\exp\int_{0}^{1}\frac{\psi(-1)}{t}dt >0
	\end{equation*}
	Thus the existence of the root $r_0$ is ensured by the Intermediate Value theorem for the continuous functions. Now if $r_0\leq 1/3$, then for the choices of functions $g=f=f_{\psi}$, equality
	$$\sum_{k=1}^{\infty}|b_k|r^k \leq d(f(0),\partial\Omega),$$
	holds; hence this prove the sharpness part. 
\end{proof}

Now we discuss a generalization of Bohr-Rogosinski sum for the class $\mathcal{F}(\psi)$.
Also see~\cite{Ganga-GBR2022,Kumar-shaoo-genBohr}. We first recall
\begin{lemma}\cite{Ganga-GBR2022}\label{series-lem-gen}
	Let $f(z)=\sum_{n=0}^{\infty}a_n z^n$ and $g(z)=\sum_{n=0}^{\infty}b_n z^n$ be analytic in $\Delta$. Let $\{\nu_{k}(r)\}_{k=0}^{\infty}$ be a sequence of non-negative functions, continuous in $[0,1)$ such that the series 
	\begin{equation*}
	\sum_{n=0}^{\infty}\vert b_n\vert \nu_{k}(r)
	\end{equation*}
	converges locally uniformly with respect to $r\in [0,1)$. If $g\prec f$, then 
	\begin{equation*}
	\sum_{n=0}^{\infty} \vert a_n\vert  \nu_n(r) \leq \sum_{n=0}^{\infty} \vert b_n\vert  \nu_n(r)
	\end{equation*}
	for all $\vert z\vert =r\leq \frac{1}{3}$.
\end{lemma}

\begin{theorem}[Generalized Bohr-Rogosinski sum]\label{Gen-Mainthm2-Psi}
	Let $\psi(\Delta)$ be convex. Let $\{\nu_{n}(r)  \}_{n=1}^{\infty}$ be a non-negative sequence of continuous functions in $[0,1]$ such that
	\begin{equation*}
	\nu_1(r)+ \sum_{n=2}^{\infty} \left\vert \frac{f^{(n)}_{\psi}(0)}{n!}\right\vert  \nu_{n}(r)
	\end{equation*}
	converges locally uniformly with respect to each $r\in[0,1)$.
	If 
	\begin{equation*}
	\vert f(z^m)\vert +\nu_{1}(r) +\sum_{n=2}^{\infty} \left\vert \frac{f^{(n)}_{\psi}(0)}{n!}\right\vert  \nu_{n}(r) < \exp \int_{0}^{1}\frac{\psi(-t)}{t}dt
	\end{equation*}
	and
	$f(z)=z+\sum_{n=2}^{\infty}a_n z^n \in \mathcal{F}(\psi)$.
	Then
	\begin{equation}\label{Gen-Mainthm2-expr}
	\vert f(z^m)\vert  + \sum_{n=1}^{\infty}\vert a_n\vert \nu_{n}(r) \leq d(0, \partial{\Omega})
	\end{equation}
	holds for $\vert z\vert =r\leq r_b=\min\{ 1/3, r_0 \}$, where $m\in \mathbb{N}$, $\Omega=f(\Delta)$ and $r_0$ is the smallest positive root of the equation:
	\begin{equation*}
	f_{\psi}(r^m) +\sum_{n=2}^{\infty} \left\vert \frac{f^{(n)}_{\psi}(0)}{n!}\right\vert  \nu_{n}(r) = \exp \int_{0}^{1}\frac{\psi(-t)}{t}dt-\nu_{1}(r),
	\end{equation*}
	where
	\begin{equation*}
	f_{\psi}(z)= z\exp\int_{0}^{z}\frac{\psi(t)}{t}dt.
	\end{equation*}
	Moreover, the inequality \eqref{Gen-Mainthm2-expr} also holds for the class $S(f)$ in $\vert z\vert =r\leq r_b$. When $r_b=r_0$, then the radius is best possible.
\end{theorem}
\begin{proof}
	Since $f\in \mathcal{F}(\psi) $, it follows from Lemma~\ref{subordination} that
	\begin{equation*}
	\frac{f(z)}{z} \prec \frac{f_{\Psi}(z)}{z}.
	\end{equation*}
	Applying Lemma~\ref{series-lem-gen}, we see that
	\begin{align}\label{Gen-Mainthm2-eq1}
	\sum_{n=1}^{\infty}|a_n|\nu_{n}(r)
	\leq \nu_{1}(r) +\sum_{n=2}^{\infty} \left|\frac{f^{(n)}_{\Psi}(0)}{n!}\right| \nu_{n}(r) \quad \text{for} \quad r\leq \frac{1}{3}.
	\end{align}
	Combining now the equation~\eqref{Gen-Mainthm2-eq1} with the Lemma~\ref{gen-thm1} gives
	\begin{align*}
	|f(z^m)|+ \sum_{n=1}^{\infty}|a_n|\nu_{n}(r)
	\leq f_{\Psi}(r^m)+ \nu_{1}(r) +\sum_{n=2}^{\infty} \left|\frac{f^{(n)}_{\Psi}(0)}{n!}\right| \nu_{n}(r)
	\end{align*}
	for $r\leq 1/3$. Hence, using the hypothesis
	\begin{equation*}
	\vert f(z^m)\vert +\nu_{1}(r) +\sum_{n=2}^{\infty} \left\vert \frac{f^{(n)}_{\psi}(0)}{n!}\right\vert  \nu_{n}(r) < \exp \int_{0}^{1}\frac{\psi(-t)}{t}dt,
	\end{equation*}
	the inequality~\eqref{Gen-Mainthm2-expr} holds in $|z|=r \leq \min \{r_0, 1/3\}$, where $r_0$ is the minimum positive root of the equation
	\begin{equation*}
	f_{\psi}(r^m) +\sum_{n=2}^{\infty} \left\vert \frac{f^{(n)}_{\psi}(0)}{n!}\right\vert  \nu_{n}(r) = \exp \int_{0}^{1}\frac{\psi(-t)}{t}dt-\nu_{1}(r),
	\end{equation*}
	In case $r_0\leq 1/3$, then for the function $f(z)=f_{\Psi}(z)=z+\sum_{n=2}^{\infty}\left|\frac{f^{(n)}_{\Psi}(0)}{n!}\right| z^n$, we have
	$$d(0, \partial{\Omega})=\exp \int_{0}^{1}\frac{\psi(-t)}{t}dt$$
	such that the equality
	\begin{equation*}
	\vert f_{\Psi}(r_{0}^{m})\vert  + \sum_{n=1}^{\infty}\vert a_n\vert \nu_{n}(r_0) =\exp \int_{0}^{1}\frac{\psi(-t)}{t}dt= d(0, \partial{\Omega})
	\end{equation*}
	holds, that is, the radius is best possible. 	
\end{proof}

In the following, taking $m\rightarrow \infty$ and $N\rightarrow 1$ gives Bohr phenomenon for the class $\mathcal{F}(\psi)$.
\begin{corollary}\label{BR-F}
	Let $\psi(\Delta)$ be convex. If
	$f(z)=z+\sum_{n=2}^{\infty}a_n z^n \in \mathcal{F}(\psi)$.
	Then
	\begin{equation}\label{Gen-Mainthm2-exprCor}
	\vert f(z^m)\vert  + \sum_{n=N}^{\infty}\vert a_n\vert |z|^n \leq d(0, \partial{\Omega})
	\end{equation}
	holds for $\vert z\vert =r\leq r_b=\min\{ 1/3, r_0 \}$, where $m\in \mathbb{N}$, $\Omega=f(\Delta)$ and $r_0$ is the smallest positive root of the equation:
	\begin{equation*}
	f_{\psi}(r^m) +\sum_{n=N}^{\infty} \left\vert \frac{f^{(n)}_{\psi}(0)}{n!}\right\vert r^n = \exp \int_{0}^{1}\frac{\psi(-t)}{t}dt,
	\end{equation*}
	where
	\begin{equation*}
	f_{\psi}(z)= z\exp\int_{0}^{z}\frac{\psi(t)}{t}dt.
	\end{equation*}
	Moreover, the inequality \eqref{Gen-Mainthm2-exprCor} also holds for the class $S(f)$ in $\vert z\vert =r\leq r_b$. When $r_b=r_0$, then the radius is best possible.
\end{corollary}

\begin{figure}[h]
	\begin{tabular}{c}
		\includegraphics[scale=0.3]{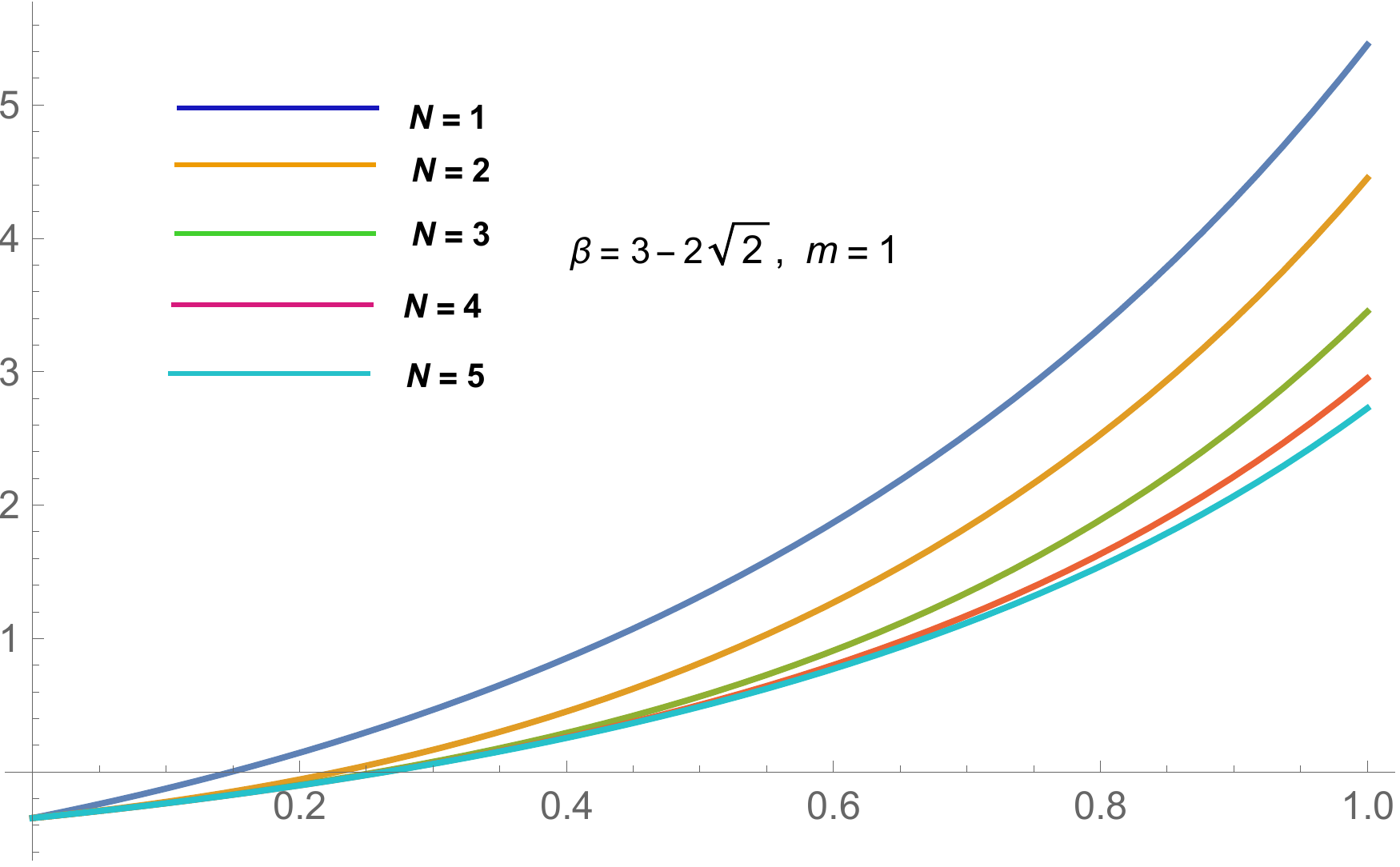}
	\end{tabular}
	\begin{tabular}{l}
		\includegraphics[scale=0.3]{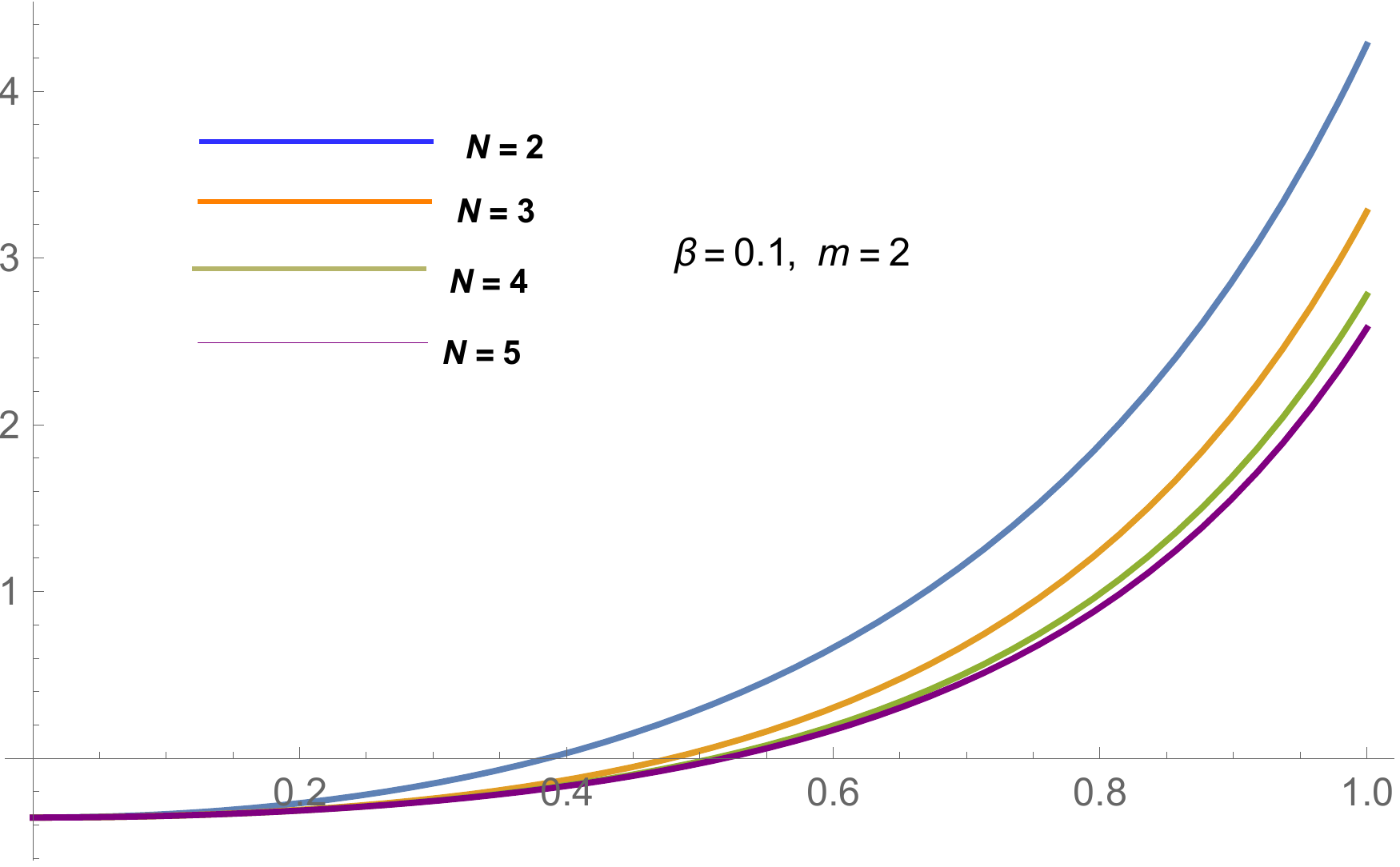}
	\end{tabular}
	\caption{Existence of roots for different values of $\beta$, $m$ and $N$ in Corollary~\ref{BR-FBooth}.} \label{BoothRogosinski}
\end{figure}

Here below, Corollary~\ref{BR-F} yields the Bohr-Rogosinki radius for the class $\mathcal{BS}(\beta)$.
\begin{corollary}\label{BR-FBooth}
	Let $0<\beta \leq 3-2\sqrt{2}$.
	If $f(z)=z+\sum_{n=2}^{\infty}a_n z^n \in \mathcal{BS}(\beta)$.
	Then
	\begin{equation}\label{Gen-Mainthm2-exprCorBooth}
	\vert f(z^m)\vert  + \sum_{n=N}^{\infty}\vert a_n\vert |z|^n \leq d(0, \partial{\Omega})
	\end{equation}
	holds for $\vert z\vert =r\leq r_b=\min\{ r_{N}^{m}(\beta), 1/3 \}$, where $m\in \mathbb{N}$, $\Omega=f(\Delta)$ and $r_{N}^{m}(\beta)$ is the smallest positive root of the equation:
	\begin{equation*}
	\hat{f}_{\psi}(r^m) -S_{N-1}(\hat{f})+ \hat{f}(r)= \left( \frac{1-\sqrt{\beta}}{1+\sqrt{\beta}} \right)^{\frac{1}{2\sqrt{\beta}}},
	\end{equation*}
	where
	\begin{equation*}
	\hat{f}(z)= z \left( \frac{1+z\sqrt{\beta}}{1-z\sqrt{\beta}} \right)^{\frac{1}{2\sqrt{\beta}}}
	\end{equation*}
	and $S_{N-1}(\hat{f})$ denotes $(N-1)^{th}$ partial sum of $\hat{f}$ and $S_0=0$.
	Moreover, the inequality \eqref{Gen-Mainthm2-exprCorBooth} also holds for the class $S(f)$ in $\vert z\vert =r\leq r_b$. The radius $r_{N}^{m}(\beta)$ is sharp with $m=1$ for all $N\in\mathbb{N}$. See, Figure~\ref{BoothRogosinski}
\end{corollary}

\begin{remark}
	Taking $m\rightarrow \infty$ and $N\rightarrow 1$ in Corollary~\ref{BR-FBooth} gives \cite[Theorem~4]{ganga-geo2021}.
\end{remark}	

We now study poly-analytic functions in connection with the class $\mathcal{F}(\psi)$. A continuous complex-valued function $F$ defined in a domain $D \subseteq \mathbb{C}$ is  poly-analytic of order $\alpha$ if it obeys the generalized Cauchy-Riemann equations $\partial^{\alpha} F/ \partial {\bar{z}}^{\alpha}=0$. An iteration argument shows that $F$ has the form
$$F(z)=\sum_{k=0}^{\alpha-1}{\bar{z}}^k f_k(z),$$
where $f_k$ are analytic functions for $k=0, \cdots, \alpha-1$. Throughout, it is assumed that $\alpha\geq2$.  In particular, a continuous complex-valued function $F=u+iv$ in a domain $D \subseteq \mathbb{C}$ is said to be bi-analytic if
$$\frac{\partial}{\partial \bar{z}} \left( \frac{\partial F}{\partial\bar{z}} \right) =0.$$
Then $F$ is a poly-analytic function of order $2$. Note that $F_{\bar{z}}$ is analytic in $D$. In any simply connected domain $D$, it can be shown that $F$ has the form
$$F(z)=\bar{z}f_1(z)+f_0(z),$$
where $f_1$ and $f_0$ are analytic functions (and obviously $f_1= F_{\bar{z}}$). Poly-analytic functions can have properties that differ from analytic functions, such as the ability to vanish on closed curves without vanishing identically. The function $F(z)=1-z\bar{z}$ is one such example.

We obtain Theorem~\ref{2} for poly-analytic function with analytic counterparts from the classes $\mathcal{F}(\psi)$, $\mathcal{S}^*(\Psi)$ and $\mathcal{C}(\Psi)$ for the Rogosinski's sum.
\begin{theorem}
	Let $F(z)=\sum_{k=0}^{\alpha-1} {\bar{z}}^k f_k(z)$ be a poly-analytic function of order $\alpha$, where $f_k$ are analytic such that $g_k=f_0 +\overline{f_k}$ is sense-preserving in $\Delta$ and $f_k(0)=0$ for each $k=0,\cdots,\alpha-1$. Suppose that
	\begin{enumerate}
		\item [$(i)$] $L(r)= B_{N}^{m}(f_\psi,r)$ if $f_0 \in \mathcal{F}(\psi)$ and $\psi(\Delta)$ is convex
		\item [$(ii)$]  $L(r)=B_{N}^{m}(f_{\Psi},r)$ if $f_0\in \mathcal{S}^{*}(\Psi)$
		\item [$(iii)$] $L(r)= B_{N}^{m}(g_{\Psi},r)$, where $zg'_{\Psi}(z)= f_{\Psi}(z)$ if $f_0 \in \mathcal{C}(\Psi)$.
		
	\end{enumerate} Then
	\begin{equation*}
	B_{N}^{m}(F,r)= \sum_{n=N}^{m} \sum_{k=0}^{\alpha-1} |a_{n,k}| r^{k+n} <1, \quad \text{if} \quad |z|<\min \{ 1/3, r_0 \}
	\end{equation*}
	where $N\in \mathbb{N}$ and $r_0$ is the minimal positive root of the equation
	\begin{equation*}
	L(r)(1-r^{\alpha})+r-1=0.
	\end{equation*}
\end{theorem}
\begin{proof}
	First see that
	\begin{align}\label{majorantsum-poly}
	B_{N}^{m}(F,r) =B_{N}^{m}\left( \sum_{k=0}^{\alpha-1}\bar{z}^k f_k(z) \right) 
	\leq \sum_{k=0}^{\alpha-1} r^k B_{N}^{m}(f_k).
	\end{align}
	Since $g_k=f_0 +\overline{f_k}$ is sense-preserving in $\Delta$ for each $k$, it follows that
	\begin{equation*}
	f'_k(z)= \omega_k(z) f'_0(z)
	\end{equation*}
	where $\omega_k$ is analytic in $\Delta$ and $|\omega_k(z)|<1$ for all $z\in \Delta$, which further using \cite[Corollary~2.1]{kamal-mediter2021} gives that for $|z|\leq 1/3$,
	\begin{align*}
	B_{N}^{m}(f_k) &= \int_{0}^{r} B_{N}^{m}(f'_k)ds
	= \int_{0}^{r} B_{N}^{m}(\omega_k f'_0)ds\\
	&\leq \int_{0}^{r}B_{N}^{m}(f'_0)ds
	=B_{N}^{m}(f_0).
	\end{align*}
	Thus, from \eqref{majorantsum-poly}
	\begin{equation*}
	B_{N}^{m}(F,r) \leq B_{N}^{m}(f_0) \left( \frac{1-r^\alpha}{1-r} \right) \quad \text{for}\quad r\leq \frac{1}{3}.
	\end{equation*}
	Now to complete the proof, we consider the following three cases:
	\begin{enumerate}
		\item [ (i) ] Let $f_0 \in \mathcal{S}^*(\psi)$. Then $f_0(z)/z \prec f_{\psi}(z)/z$ holds, which with \cite[Lemma~2.1]{kamal-mediter2021} gives
		\begin{equation*}
		B_{N}^{m}(f_0) \leq B_{N}^{m}(f_{\psi}) \quad \text{for} \quad r\leq \frac{1}{3}.
		\end{equation*}
		\item [ (ii) ]  Let $f_0 \in \mathcal{C}(\psi)$. Then using Alexander's relation, $f_\psi(z)=zg'_{\psi}(z)$ it follows from first case that
		\begin{equation*}
		B_{N}^{m}(f_0) \leq B_{N}^{m}(g_{\psi}) \quad \text{for} \quad r\leq \frac{1}{3}.
		\end{equation*}
		\item [ (iii) ] Let $\Psi$ be convex and $f_0 \in \mathcal{F}(\Psi)$. Then from Lemma~\ref{subordination} , it follows that $f_0(z)/z \prec f_{\Psi}(z)/z$. This using \cite[Lemma~2.1]{kamal-mediter2021} gives that
		\begin{equation*}
		B_{N}^{m}(f_0) \leq B_{N}^{m}(f_{\Psi}) \quad \text{for} \quad r\leq \frac{1}{3}.
		\end{equation*}
		Hence, it follows that
		\begin{align*}
		B_{N}^{m}(F,r) \leq L(r) \left(\frac{ 1-r^{\alpha}}{1-r} \right) 
		< 1
		\end{align*}
		for $|z|=r < \min\{1/3, r_0 \}$, where $r_0$ is the minimum positive root of the polynomial equation
		\begin{equation*}
		L(r)(1-r^{\alpha})+r-1=0,
		\end{equation*}
		and $L(r)$ is as defined in the hypothesis. 
	\end{enumerate}
\end{proof}	

\begin{corollary}
	Let $F(z)=\sum_{k=0}^{\alpha-1} {\bar{z}}^k f_k(z)$ be a poly-analytic function of order $\alpha$, where $f_k$ are analytic such that $g_k=f_0 +\overline{f_k}$ is sense-preserving in $\Delta$ and $f_k(0)=0$ for each $k=0,\cdots,\alpha-1$. Let $\beta \in (0,3-2\sqrt{2}]$ and suppose that $f_0 \in \mathcal{BS}(\beta)$. Then
	\begin{equation*}
	\sum_{n=1}^{m} \sum_{k=0}^{\alpha-1} |a_{n,k}| r^{k+n} <1, \quad \text{if} \quad |z|< 1/3.
	\end{equation*}
\end{corollary}

\begin{corollary}
	Let $F(z)=\sum_{k=0}^{\alpha-1} {\bar{z}}^k f_k(z)$ be a poly-analytic function of order $\alpha$, where $f_k$ are analytic such that $g_k=f_0 +\overline{f_k}$ is sense-preserving in $\Delta$ and $f_k(0)=0$ for each $k=0,\cdots,\alpha-1$. Then
	\begin{equation*}
	B_{N}^{m}(F,r)= \sum_{n=N}^{m} \sum_{k=0}^{\alpha-1} |a_{n,k}| r^{k+n} <1, \quad \text{if} \quad |z|< \min\{ 1/3, r_0 \}
	\end{equation*}
	where $N\in \mathbb{N}$ and $r_0$ is the minimal positive root of the equations:
	\begin{equation*}
	\left (\frac{r}{1-r} \right)(r^N -r^m)(1-r^{\alpha})+r-1=0,  \quad \text{if} \quad  f_0 \in \mathcal{C}.
	\end{equation*}
	and
	\begin{equation*}
	\left(\frac{r}{1-r} \right)\left[\frac{r^{N}-r^{m}}{r(1-r)} -(m-1)r^m +(N-1)r^N \right] (1-r^{\alpha})+r-1=0, \; \text{if} \quad f_0 \in \mathcal{S}^{*}.
	\end{equation*}
\end{corollary}		

\begin{remark}\label{rogosinski-convex}
	Let $F(z)=\sum_{k=0}^{\alpha-1} {\bar{z}}^k f_k(z)$ be a poly-analytic function of order $\alpha$. Assume that $f_k \in \mathcal{C}$  for each $k=0,\cdots,\alpha-1$. Then the Rogosinski sum
	\begin{equation*}
	B_{1}^{m}(F,r)= \sum_{n=1}^{m} \sum_{k=0}^{\alpha-1} |a_{n,k}| r^{k+n} <1 \quad \text{in} \quad |z|< R_{\alpha}
	\end{equation*}
	where the Rogosinski radius $R_{\alpha}$ is the minimal positive root of the equations
	\begin{equation*}
	r(1-r^{\alpha})-(1-r)^2=0.
	\end{equation*}
	Here, we note that the radius $R_{\alpha}>1/3$.
\end{remark}

Motivated from the Remark~\ref{rogosinski-convex}, we discuss the result when $f_k$ comes from the Janowski-starlike class.
\begin{theorem}\label{Poly-Janowski}
	Let $F(z)=\sum_{k=0}^{\alpha-1} {\bar{z}}^k f_k(z)$ be a poly-analytic function of order $\alpha$. Assume that $f_k \in \mathcal{S}^*[D,E]$  for each $k=0,\cdots,\alpha-1$. Then the Rogosinski sum
	\begin{equation*}
	B_{1}^{m}(F,r)= \sum_{n=1}^{m} \sum_{k=0}^{\alpha-1} |a_{n,k}| r^{k+n} <1 \quad \text{in} \quad |z|< R
	\end{equation*}
	for each $m\geq 2$, where the Rogosinski radius $R$ is the minimal positive root of the equations
	\begin{equation*}
	f_j(r)(1-r^{\alpha})+r-1=0,
	\end{equation*}
	where
	\begin{equation}\label{Janow-extFunct}
	f_j(z)= 
	\left\{
	\begin{array}
	{lr}
	z(1+E z)^{\frac{D-E}{E}}, & E\neq0; \\
	ze^{Dz},   & E=0.
	\end{array}
	\right.
	\end{equation}
\end{theorem}
\begin{proof}
	Since for each $k$, it is well-known that
	\begin{equation*}
	|a_{n,k}| \leq \prod_{t=0}^{n-2}\frac{|E-D+Et|}{t+1},
	\end{equation*}
	where the equality is achieved for the function $f_j$ given by \eqref{Janow-extFunct}.
	Now 
	\begin{align*}
	B_{N}^{m}(F,r) &=\sum_{k=0}^{\alpha-1} r^k \left(\sum_{n=N}^{m}|a_{n,k}|r^n \right)\\
	& \leq \sum_{k=0}^{\alpha-1} r^k \left(\sum_{n=N}^{m} \prod_{t=0}^{n-2}\frac{|E-D+Et|}{t+1} r^n \right)\\
	&=\sum_{k=0}^{\alpha-1} r^k  B_{N}^{m}(f_j,r)
	= B_{N}^{m}(f_j) \sum_{k=0}^{\alpha-1} r^k
	= B_{N}^{m}(f_j) \left( \frac{1-r^{\alpha}}{1-r} \right).
	\end{align*}
	Therefore, 
	$B_{N}^{m}(F,r) <1$
	holds in $|z|<R_{N}^{m}$, where $R_{N}^{m}$ is the minimum positive root of the polynomial equation
	\begin{equation*}
	B_{N}^{m}(f_j,r)(1-r^{\alpha})+r-1=0.
	\end{equation*}
	Since $B_{N}^{m}(f_j) \leq B_{N}^{m+1}(f_j)$. Therefore, taking $N=1$ and $m \rightarrow \infty$ shows that $R_{N}^{m}$ approaches the Rogosisnki radius $R$. 
\end{proof}

\begin{corollary}\label{particular-starlikeconvex}
	Let $F(z)=\sum_{k=0}^{\alpha-1} {\bar{z}}^k f_k(z)$ be a poly-analytic function of order $\alpha$, where $f_k$ are analytic such that $g_k=f_0 +\overline{f_k}$ is sense-preserving in $\Delta$ and $f_k(0)=0$ for each $k=0,\cdots,\alpha-1$. Then
	\begin{equation*}
	B_{1}^{m}(F,r)= \sum_{n=1}^{m} \sum_{k=0}^{\alpha-1} |a_{n,k}| r^{k+n} <1, \quad \text{if} \quad |z|< \min\{ 1/3, r_0 \}
	\end{equation*}
	where $r_0$ is the minimal positive root of the equations:
	\begin{equation*}
	\left (\frac{r}{1-r} \right)(r -r^m)(1-r^{\alpha})+r-1=0,  \quad \text{if} \quad  f_0 \in \mathcal{C}.
	\end{equation*}
	and
	\begin{equation*}
	\left(\frac{r}{1-r} \right)\left[\frac{r-r^{m}}{r(1-r)} -(m-1)r^m \right] (1-r^{\alpha})+r-1=0, \quad \text{if} \quad  f_0 \in \mathcal{S}^{*}.
	\end{equation*}
\end{corollary}		
\begin{remark}
	In Corollary~\ref{particular-starlikeconvex}, when $f_0 \in \mathcal{C}$. Then for each value of $m, \alpha \geq2$, 
	\begin{equation*}
	\sum_{n=1}^{m} \sum_{k=0}^{\alpha-1} |a_{n,k}| r^{k+n} <1 \quad \text{in} \quad |z|< \frac{1}{3}.
	\end{equation*}
	When $f_0 \in \mathcal{S}^*$, then $r_0 \leq 1/3.$
\end{remark}

\begin{corollary}
	Let $F(z)=\sum_{k=0}^{\alpha-1} {\bar{z}}^k f_k(z)$ be a poly-analytic function of order $\alpha$, where $f_k$ are analytic such that $g_k=f_0 +\overline{f_k}$ is sense-preserving in $\Delta$ and $f_k(0)=0$ for each $k=0,\cdots,\alpha-1$. Then
	\begin{equation*}
	B_{1}^{m}(F,r)= \sum_{n=1}^{\infty} \sum_{k=0}^{\alpha-1} |a_{n,k}| r^{k+n} <1, \quad \text{if} \quad |z|< \min\{1/3,r_0\}
	\end{equation*}
	where $r_0$ is the minimal positive root of the equations:
	\begin{equation*}
	(1-r)^2-r^2+r^{\alpha+2}=0,  \quad \text{if} \quad  f_0 \in \mathcal{C}
	\end{equation*}
	and
	\begin{equation*}
	(1-r)^3 -r+r^{\alpha+1}=0, \quad \text{if} \quad  f_0 \in \mathcal{S}^{*}.
	\end{equation*}
\end{corollary}		

We now conclude this section with an interesting class studied by Bhowmik~\cite{Bhowmik-Concave2012}  on the basis of geometry, also see the family of Generalized Pascal Snails~\cite{Kanas-pascal}, which is as follows:
\begin{definition}\cite{Bhowmik-Concave2012}
	A function $f: \Delta \rightarrow \mathbb{C}$ is said to be in the class of concave univalent functions with opening angle $\pi \beta$, $\beta \in [1,2]$, at infinity if it fulfills the following conditions:
	\begin{enumerate}
		\item [$(a)$] $f$ is analytic in $\Delta$ with the standard normalization and $f(1)= \infty$.
		\item  [$(b)$] $f$ maps $\Delta$ conformally onto a set whose complement with respect to $\mathbb{C}$ is convex.
		\item  [$(c)$] the opening angle of $f(\Delta)$ at $\infty$ is less than or equal to $\pi \beta$.
	\end{enumerate}
\end{definition}

For the function $f\in Co(\beta)$, the image domain $\mathbb{C}/f(\Delta)$ is closed and unbounded. Bhowmik~\cite{Bhowmik-Concave2012} obtained sharp Coefficient's bounds, growth theorem, rotation theorem and distortion theorem etc. with the extremal function given by
\begin{equation}\label{cocave-extremal}
f_{\beta}(z):=\frac{1}{2\beta}\bigg[\left(\frac{1+z}{1-z} \right)^{\beta} -1 \bigg].
\end{equation}
It seems interesting to study as an independent interest the class defined by
\begin{equation*}
\mathcal{F}(f_{\beta}) := \left\{ f\in \mathcal{A} : \left(\frac{zf'(z)}{f(z)}-1 \right) \prec f_{\beta}(z) \right\}.
\end{equation*}
For the class $Co(\beta)$, Bhowmik and Das~\cite{bhowmik2018} studied the Bohr phenomenon. Now we explore this in connection with poly-analytic functions. 
\begin{theorem}\label{Poly-concave}
	Let $F(z)=\sum_{k=0}^{\alpha-1} {\bar{z}}^k f_k(z)$ be a poly-analytic function of order $\alpha$, where $f_k\in Co(\beta) $. Then
	\begin{equation*}
	B_{N}^{m}(F,r)= \sum_{n=N}^{m} \sum_{k=0}^{\alpha-1} |a_{n,k}| r^{k+n} <1, \quad \text{if} \quad |z|< r_{\beta}(N)
	\end{equation*}
	where $N\in \mathbb{N}$ and $r_{\beta}(N)$ is the minimal positive root of the equation
	\begin{equation*}
	(1-r^{\alpha}) B_{N}^{m}(f_{\beta},r)+r-1=0,
	\end{equation*}
	where $f_{\beta}$ is given by $\eqref{cocave-extremal}$.
\end{theorem}
\begin{proof}
	First note that
	\begin{align*}
	B_{N}^{m}(F,r) &=\sum_{k=0}^{\alpha-1} r^k \left(\sum_{n=N}^{m}|a_{n,k}|r^n \right)
	\leq \sum_{k=0}^{\alpha-1} r^k \left(\sum_{n=N}^{m} A_n r^n \right)=\sum_{k=0}^{\alpha-1} r^k  B_{N}^{m}(f_{\beta},r) \\
	&= B_{N}^{m}(f_{\beta}) \sum_{k=0}^{\alpha-1} r^k 
	= B_{N}^{m}(f_{\beta}) \left( \frac{1-r^{\alpha}}{1-r} \right).
	\end{align*}
	Therefore, 
	\begin{equation*}
	B_{N}^{m}(F,r) <1
	\end{equation*}
	holds in $|z|<r_{\beta}(N)$, where $r_{\beta}(N)$ is the minimum positive root of the polynomial equation
	\begin{equation*}
	(1-r^{\alpha}) B_{N}^{m}(f_{\beta},r)+r-1=0.
	\end{equation*}
	This completes the proof. 
\end{proof}	
\begin{remark}
	Taking $\beta=1$, in Theorem~\ref{Poly-concave} yields the special case for convex functions.
\end{remark}

From the Proof of Theorem~\ref{Poly-concave}, we deduce that
\begin{corollary}
	Let $F(z)=\sum_{k=0}^{\alpha-1} {\bar{z}}^k f_k(z)$ be a poly-analytic function of order $\alpha$, where $f_k \in \mathcal{F}(f_{1}) $. Then
	\begin{equation*}
	|F(z)|+ \sum_{n=N}^{\infty} \sum_{k=0}^{\alpha-1} |a_{n,k}| r^{k+n} <1, \quad \text{if} \quad |z|< \min\{r_0, 1/3 \}
	\end{equation*}
	where $N\in \mathbb{N}$ and $r_0$ is the minimal positive root of the equation
	\begin{equation*}
	(1-r^{\alpha})(r+r^N)-(1-r)^2=0.
	\end{equation*}
	As $\alpha\rightarrow \infty$ and $N=2$, then $r_0\rightarrow 1/3$ and is best possible in this case.
\end{corollary}			

\section{ Classical coefficient-functionals bounds for $\mathcal{F}(\psi)$}
Suppose that the function $\psi$ is univalent in $\Delta$ with $\psi(0)=0$. Then we can write
\begin{equation}\label{powerseries-Univalent}
\psi(z)=A_1 z+A_2z^2+A_3z^3+ \cdots.
\end{equation}
Further, if the function $1+\psi(z)$ belongs to the Carath\'{e}odory class then the Taylor's coefficients satisfies the well known bounds $|A_n|\leq 2$. Otherwise, when $A_1\neq 0$ then we have $|A_n|/|A_1| \leq n$ for each $n\in \mathbb{N}$. Since the Taylor's coefficients of the function $f(z)=z+\sum_{k=2}^{\infty}a_k z^k$ in the class $\mathcal{F}(\psi)$ can be expressed in the form of $A_i$, therefore to study the behavior of $|a_k|$ and some other coefficient functional, it is natural to consider the normalized form of $f$, that is $A_1=1$. 

We start with some basic transformations which hold also for the class $\mathcal{S}$ of normalized univalent functions.
\begin{remark}
	Let $f \in \mathcal{F}(\psi)$, where $\psi\in \mathcal{S}$. Then the following transformations belong to $\mathcal{F}(\psi)$:
	\begin{enumerate}[$(i)$]
		\item Let $\alpha \in \mathbb{R}$,
		\begin{equation*} 
		e^{i \alpha} f(e^{i \alpha} z)= z+ \sum_{n=2}^{\infty} a_n e^{i(n-1)\alpha} z^n.
		\end{equation*}
		\item Let $0<t\leq 1$,
		\begin{equation*} 
		\frac{1}{t} f(t z)=z+ \sum_{n=2}^{\infty} a_n t^{n-1} z^n.
		\end{equation*}
		
		\item 	Let $k$ is an positive integer
		\begin{equation*} 
		[f(z^k)]^{\frac{1}{k}}= z+ \frac{a_2}{k} z^{k+1}+ \frac{1}{2k^2}(2ka_3-(k-1)a_{2}^{2}) z^{2k+1} + \cdots.
		\end{equation*}
	\end{enumerate}
\end{remark}

First, we obtain sharp initial coefficient's bounds, Zalcman~\cite{ChoKwonLecko-Zalc2018,Ma-Zalcman1999} and  the Fekete-Szeg\"{o}~\cite{AliRavi-Pvalent2007,minda94} functional.
\begin{theorem}
	If $f \in \mathcal{F}(\psi)$, where $\psi\in \mathcal{S}$. Then 
	\begin{enumerate}
		\item	sharp bounds for initial coefficient's are
		\begin{equation*}
		|a_2|\leq A_1, \quad	 	
		2|a_3| \leq \left\{
		\begin{array}
		{lr}
		{A_2+A^2_1}, & \text{if} \quad {A^2_1+A_2} \geq A_1\\
		A_1, & \text{if} \quad A_1 \geq A^2_1+A_2 \geq -A_1\\
		-(A_2+A^2_1), & \text{if} \quad {A^2_1+A_2} \leq -A_1
		\end{array}
		\right.	 	
		\end{equation*}
		and
		\begin{equation*}
		|a_4| \leq \frac{A_1}{3}{H(q_1, q_2)}, \quad\text{where} \quad 
		\left\{
		\begin{array}
		{lr}
		q_1= \frac{4A_2+3A^2_1}{2A_1} &  \\
		q_2= \frac{2A_3+A^3_1+3A_1A_2}{2A_1} & 
		\end{array}
		\right.	 	
		\end{equation*}
		\item	sharp bound for the Zalcman functional of early coefficients is
		\begin{equation*}
		|a_2 a_3 - a_4| \leq \frac{A_1}{3}{H(q_1, q_2)}, \quad\text{where} \quad 
		\left\{
		\begin{array}
		{lr}
		q_1= \frac{A^2_1 +4A_2}{2A_1} &  \\
		q_2= \frac{A_3-A^2_1}{A_1} & 
		\end{array}
		\right.	 	
		\end{equation*}
		\item	 sharp bound for the Fekete-Szeg\"{o} functional is
		\begin{equation*}
		|a_3- \nu a_{2}^{2}| \leq \left\{
		\begin{array}
		{lr}
		\frac{1}{2}(A_2 + (1-2\nu)A^2_1), & \text{if} \quad A_1 + (2\nu -1)A^2_1 \leq A_2\\
		A_1/2, & \text{if} \quad  \frac{A_2-A_1+A^2_1}{2A^2_1} \leq \nu \leq \frac{A_2 +A_1+A^2_1}{2A^2_1}\\
		-\frac{1}{2}(A_2 + (1-2\nu)A^2_1), & \text{if} \quad (2\nu-1)A^2_1 \geq A_1+A_2
		\end{array}
		\right.	 	
		\end{equation*}  
		
	\end{enumerate}	
	where  $H(q_1,q_2)$ and $(q_1,q_2)$ are given by \cite[Lemma~2]{Prokhorov-1981}.
\end{theorem}
\begin{proof}
	The proof follows by employing the technique of \cite[Theorem~1]{AliRavi-Pvalent2007}, a lemma of Ma and Minda~\cite{minda94} and \cite[Lemma~2]{Prokhorov-1981}. Here in this proof, we simply notice that $|A_3|\leq 3$, which is not the case for $\Psi$ in the class $\mathcal{S}^*(\Psi)$. 	
\end{proof}	

\begin{corollary}\label{BoothCor}
	Let $f \in \mathcal{BS}(\beta) $, $0\leq \beta <1$. Then
	\begin{enumerate}
		\item	sharp bounds for initial coefficient's are
		\begin{equation*}
		|a_2|\leq 1, \quad	 	
		|a_3| \leq 1/2
		\end{equation*}
		and
		\begin{equation*}
		|a_4| \leq \frac{1}{3} H\left(\frac{3}{2}, \frac{2\beta+1}{2} \right)
		\end{equation*}
		\item	sharp bound for the Zalcman functional of early coefficients is
		\begin{equation*}
		|a_2 a_3 - a_4| \leq \frac{1}{3}{H(1/2, \beta-1)}
		\end{equation*}
		\item	 sharp bound for the Fekete-Szeg\"{o} functional is
		\begin{equation*}
		2|a_3- \nu a_{2}^{2}| \leq \left\{
		\begin{array}
		{lr}
		(1-2\nu), & \text{if} \quad \nu \leq 0\\
		1, & \text{if} \quad  0\leq \nu \leq 1\\
		(2\nu-1), & \text{if} \quad \nu \geq 1
		\end{array}
		\right.	 	
		\end{equation*}  
		
	\end{enumerate}	
	where  $H(q_1,q_2)$ is given by \cite[Lemma~2]{Prokhorov-1981}. The equality cases follow for the extremal function
	\begin{equation*}
	f_{\beta}(z)= \left( \frac{1+z \sqrt{\beta}}{1-z \sqrt{\beta}} \right)^{\frac{1}{2\sqrt{\beta}}} = z+z^2+ \frac{1}{2}z^3+ \frac{1}{3}\left( \beta+\frac{1}{2} \right) z^4 +\cdots.
	\end{equation*}
\end{corollary}	
\begin{remark}
	In Corollary~\ref{BoothCor}, part-$1$ improves $|a_3|$ and $|a_4|$ bounds in \cite[Theorem~2.2]{kargar-2019}, while part-$3$ reduces to \cite[Corollary~2.2]{kargar-2019}.
\end{remark}

\begin{corollary}
	Let $f \in \mathcal{S}_{cs}(\beta)$, $0 \leq \beta <1$. Then
	\begin{enumerate}
		\item	sharp bounds for initial coefficient's are
		\begin{equation*}
		|a_2|\leq 1, \quad	 	
		2|a_3| \leq 2-\beta
		\end{equation*}
		and
		\begin{equation*}
		|a_4| \leq \frac{1}{3} H \left( \frac{7-4\beta}{2}, \frac{1+{\beta}^2}{1+\beta}+\frac{4-3\beta}{2} \right)
		\end{equation*}
		\item	sharp bound for the Zalcman functional of early coefficients is
		\begin{equation*}
		|a_2 a_3 - a_4| \leq \frac{1}{3} H \left( \frac{5-4\beta}{2}, \beta(\beta-1) \right)
		\end{equation*}
		\item	 sharp bound for the Fekete-Szeg\"{o} functional is
		\begin{equation*}
		2|a_3- \nu a_{2}^{2}| \leq \left\{
		\begin{array}
		{lr}
		2(1-\nu)-\beta, & \text{if} \quad 1-2\nu \geq \beta\\
		1, & \text{if} \quad  1-\beta \leq 2\nu \leq 3-\beta\\
		\beta-2(1-\nu), & \text{if} \quad 2\nu \geq 3-\beta
		\end{array}
		\right.	 	
		\end{equation*}  
		
	\end{enumerate}	
	where  $H(q_1,q_2)$ is given by \cite[Lemma~2]{Prokhorov-1981}. The equality cases follow for the extremal function
	\begin{equation*}
	f_{\beta}(z)=z\left( \frac{1+\beta z}{1-z} \right)^{\frac{1}{1+\beta}}= z+z^2+ \left( 1-\frac{\beta}{2} \right) z^3+ \frac{1}{6}(6-5\beta+2 {\beta}^2) z^4 + \cdots.
	\end{equation*}
\end{corollary}	

Using the coefficients of the analytic function $f$ given by  $f(z)=z+\sum_{k=n+1}^{\infty}b_{k}z^{k}$, Pommerenke, and Noonan and Thomas considered the Hankel determinant $H_q(n)$ which is defined by
\begin{equation*}
H_q(n):=\left|\begin{matrix}
b_n     & b_{n+1} & \dots & b_{n+q-1}\\
b_{n+1} & b_{n+2} & \dots & b_{n+q}\\
\vdots  & \vdots  & \ddots & \vdots\\
b_{n+q-1}&b_{n+q} & \dots & b_{n+2(q-1)}
\end{matrix}\right|,
\end{equation*}
where $b_1=1$. Estimation of upper bound of $|H_3(1)|$~\cite{Leko8/9,Lecko1/2}, $|H_2(2)|$  and $|H_2(1)|$~\cite{AliRavi-Pvalent2007} for the functions belonging to various subclasses of univalent functions is a usual phenomenon in GFT, see~\cite{Kumar-cardioid}. For results on $\mathcal{S}^*(\Psi)$, we refer to see~\cite{AliRavi-Pvalent2007,Alarif-H222017,Lee-secHankl2013,Ma-Zalcman1999}

We next tackle second hankel determinant $H_2(2)= (a_2a_4-a^2_3 )$ for the class $\mathcal{F}(\psi)$: Here one should note that we have $|A_3|\leq 3$, which is not the case in Ma-Minda function $\Psi$.
\begin{theorem}\label{secondhankel}
	If $f \in \mathcal{F}(\psi)$, where $\psi\in \mathcal{S}$.
	\begin{enumerate}
		\item [$(i)$] If $|A_2|\leq A_1$ and $|A^4_1-4A_1A_3+3A^2_2+6A^2_1A_2|\leq 3 A^2_1$ holds , then 
		$$|H_2(2)| \leq A^2_1/4.$$
		
		\item [$(ii)$] If either $$|A_2|\geq A_1 \quad \text{and} \quad |A^4_1-4A_1A_3+3A^2_2+6A^2_1A_2|\geq A_1|A_2|+ 2A^2_1,$$
		or the conditions
		$$|A_2|\leq A_1 \quad \text{and} \quad  |A^4_1-4A_1A_3+3A^2_2+6A^2_1A_2| \geq 3 A^2_1 $$ 
		holds, then 
		\begin{equation*}
		|H_2(2)| \leq \frac{|A^4_1-4A_1A_3+3A^2_2+6A^2_1A_2| }{12}
		\end{equation*}
		
		\item [$(iii)$] If $|A_2| > A_1$ and $2A_1|A_2|+2A^2_1 \geq |A^4_1-4A_1A_3+3A^2_2+6A^2_1A_2|$ holds, then
		\begin{equation*}
		|H_2(2)| \leq \frac{A^2_1}{4}-\frac{4A^2_1(|A_2|-A_1)^2}{48(|A^4_1-4A_1A_3+3A^2_2+6A^2_1A_2|-2A_1|A_2|-A^2_1)}.
		\end{equation*}
	\end{enumerate}
\end{theorem}
\begin{proof}
	Let us consider the following series expansion:
	\begin{equation*}
	\frac{zf'(z)}{f(z)}= 1+ a_2z +(2a_3-a^2_2)z^2 + \cdots.
	\end{equation*}
	Since
	\begin{equation*}
	\frac{zf'(z)}{f(z)}-1=\psi(\omega(z))=\psi\left( \frac{p(z)+1}{p(z)-1} \right).
	\end{equation*}
	Therefore, we get
	\begin{equation*}
	2a_2={A_1p_1}, \quad 8a_3={A_2p^2_1 -A_1(p^2_1-2p_2) +A^2_1p^2_1}
	\end{equation*}
	and
	\begin{align*}
	a_4=&\frac{1}{48}\bigg( 2p_1(-2A_2p^2_1 +A_3p^2_1+4A_2p_2) +3A^3_1p^3_1 -A^2_1p_1(p^2_1-2p_2)\\
	&+A_1\bigg(p^3_1(6+2(A_2-1)) +3(A_2-1)+1-2A_2 \bigg) -8p_1p_2 +8p_3 \bigg)
	\end{align*}
	Now using the expressions from Lemma for $p_2$ and $p_3$ gives
	\begin{align}\label{secondhankel-expression}
	a_2a_4-a^2_3 = &\frac{1}{192} \bigg( M p^4  +2p^2 x(4-p^2)A_1 A_2 -4p^2x^2(4-p^2)A^2_1 \nonumber\\
	&-3(4-p^2)^2 x^2 A^2_1 +8p(4-p^2)z(1-|x|^2)A^2_1 \bigg),
	\end{align}
	where $M:=-A^4_1+4A_1A_3-3A^2_2-6A^2_1A_2$.  Assume $p_1=p\in[0,2]$ because the function $P(e^{i\theta}z)$ is Carath\'{e}odory for each $0\leq \theta<2\pi$ if $P(z)$ is so. Now applying the triangle inequality in \eqref{secondhankel-expression} with $|x|=\rho$ gives
	\begin{align*}
	|a_2a_4-a^2_3| &\leq \frac{1}{192} \bigg( |M|p^4 +4p^2(4-p^2){\rho}^2 A^2_1 +8A^2_1p(4-p^2)(1-{\rho}^2)  \\
	&\quad +2A_1|A_2|p^2 (4-p^2)\rho +3A^2_1(4-p^2)^2{\rho}^2 \bigg)\\
	&=: G(p,\rho).
	\end{align*}
	It is an easy exercise to note that $\max_{0\leq\rho\leq1}G(p,\rho)=G(p,1)=:F(p)$. Further, a simplication gives
	\begin{align*}
	192F(p)
	= & \left(|M| -A_1(3xA_1 -2|A_2|-4A_1) \right)p^4 +8A_1\left(|A_2|-A_1 \right)p^2 +16A^2_1\\
	=: &\left( Ap^4 +Bp^2+C \right).
	\end{align*}
	The result follows by applying the known inequality
	\begin{align*}
	\max_{0\leq t\leq4}\left( At^4 +Bt^2+C \right)
	= 
	\left\{
	\begin{array}
	{lr}
	C, & \quad \text{if} \quad B\leq0, 4A\leq{-B} \\
	16A+4B+C,  &  \quad \text{if} \quad B\geq0, 8A\geq {-B}, \quad\text{or}\quad B\leq 0, 4A\geq {-B}\\
	\frac{4AC-B^2}{4A}, & \quad \text{if} \quad B>0, 8A\leq {-B}.
	\end{array}
	\right.
	\end{align*}
	This completes the proof. 
\end{proof}
\begin{remark}
	It is worth here to mention that the conclusion of the Theorem~\ref{secondhankel-expression} also holds for the class $\mathcal{S}^*(\Psi)$ of Ma-Minda starlike function when $\Psi(z):=1+\psi(z)$ is a Carath\'{e}odory function and $\psi$ has the form \eqref{powerseries-Univalent}.
\end{remark}

\begin{corollary}
	Let $f\in \mathcal{BS}(\beta)$, $0 \leq \beta <1$. Then for all $\beta$
	\begin{equation*}
	|H_2(2)|= |a_2a_4-a^2_3| \leq \frac{1}{4}.
	\end{equation*}
	The inequality is sharp with the extremal function
	\begin{equation*}
	f(z)= z\exp\left( \int_{0}^{z}\frac{t}{1-\beta t^4} dt \right) =z+ \frac{z^3}{2}+ \frac{z^5}{8} + \left( \frac{1+8 \beta}{48} \right) z^7 +\cdots.
	\end{equation*}
\end{corollary}			

\begin{corollary}
	Let $f \in \mathcal{S}_{cs}(\beta)$, $0 \leq \beta <1$. Then
	\begin{equation*}
	|H_2(2)|= |a_2a_4-a^2_3| \leq \frac{1}{4}, \quad\text{if} \quad |6-8\beta-{\beta}^2| \leq 3
	\end{equation*}
	and 
	\begin{equation*}
	|H_2(2)| \leq \frac{|6-8\beta-{\beta}^2|}{12}, \quad\text{if} \quad |6-8\beta-{\beta}^2| \geq 3.
	\end{equation*}
	The equality in first case holds with the extremal function
	\begin{equation*}
	f(z)= z\exp\left( \int_{0}^{z}\frac{t}{1+(\beta-1)t^2-\beta t^4} dt \right) =z+ \frac{z^3}{2}+ \left( \frac{3- 2\beta}{8} \right) z^5 +\cdots.
	\end{equation*}
\end{corollary}		

Further, application of this section can be seen for the non-Carath\'{e}odory regions given by Generalized Pascal Snails~\cite{Kanas-pascal}.

\end{document}